\documentclass[12pt]{article}
\usepackage{amssymb}
\usepackage{graphicx}
\usepackage{amsmath}

\newtheorem{theorem}{Theorem}

\newtheorem{corollary}[theorem]{Corollary}

\newtheorem{lemma}[theorem]{Lemma}

\newenvironment{proof}[1][Proof]{\textbf{#1.} }{\ \rule{0.5em}{0.5em}}

\begin{document}

\title{On the probability of satisfying a word in a group}
\author{Mikl\'{o}s Ab\'{e}rt}
\maketitle

\begin{abstract}
We show that for any finite group $G$ and for any $d$ there exists a word $%
w\in F_{d}$ such that a $d$-tuple in $G$ satisfies $w$ if and only if it
generates a solvable subgroup. In particular, if $G$ itself is not solvable,
then it cannot be obtained as a quotient of the one relator group $%
F_{d}/\left\langle w\right\rangle $.

As a corollary, the probability that a word is satisfied in a fixed
non-solvable group can be made arbitrarily small, answering a question of
Alon Amit.
\end{abstract}

\section{Introduction}

Let $F_{n}$ denote the free group on $n$ letters and let $G$ be a group. For 
$w\in F_{n}$ we say that the $n$-tuple $(g_{1},g_{2},\ldots ,g_{n})\in G^{n}$
\emph{satisfies} $w$ if the substitution $w(g_{1},g_{2},\ldots ,g_{n})=1$.
Our first result is the following.

\begin{theorem}
\label{nagy}Let $G$ be a finite group. Then for all $n$ there exists a word $%
w\in F_{n}$ such that for all $g_{1},g_{2},\ldots ,g_{n}\in G$, the tuple $%
(g_{1},g_{2},\ldots ,g_{n})$ satisfies $w$ if and only if the subgroup $%
\left\langle g_{1},g_{2},\ldots ,g_{n}\right\rangle \leq G$ is solvable.
\end{theorem}

Note that if $G$ itself is not solvable, then a word as in Theorem \ref{nagy}
has to use at least $n-2$ letters of $F_{n}$. Indeed, if $w$ omits at least
two letters, then any two elements of $G$ generates a solvable subgroup,
which, using a theorem of Thompson \cite{thom} (see also \cite{flav})
implies that $G$ itself is solvable.

For $w\in F_{n}$ let $F_{n}/\left\langle w\right\rangle $ denote the
one-relator group defined by $w$. As an immediate corollary of Theorem \ref
{nagy}, we get the following.

\begin{corollary}
Let $G$ be a finite non-solvable group. Then for all $n$ there exists a word 
$w\in F_{n}$ such that $G$ is not a quotient of $F_{n}/\left\langle
w\right\rangle $.
\end{corollary}

We suspect that this property holds exactly when $G$ is not solvable.

\bigskip

\noindent \textbf{Question 1. }\textit{Let }$G$\textit{\ be a finite
solvable group. Does there exist }$k\in \mathbb{N}$\textit{\ such that for
all }$n\geq k$\textit{\ and every }$w\in F_{n}$\textit{\ the one-relator
group }$F_{n}/\left\langle w\right\rangle $\textit{\ has a surjective
homomorphism to }$G$\textit{?}

\bigskip

For $w\in F_{n}$ let $P(G,w)$ denote the probability that for $n$
independent uniform random elements $g_{1},\ldots ,g_{n}\in G$ we have $%
w(g_{1},\ldots ,g_{n})=1$. Note that $P(G,w)$ only depends on the word $w$
and not on $n$ so we can assume $w\in F_{\infty }$.

The probabilities $P(G,w)$ have been investigated in the literature mainly
for a fixed word and running $G$. The strongest result in this direction is
of Dixon, Pyber, Seress and Shalev \cite{dixetal} who proved that for any
fixed word $1\neq w$ the probability $P(G,w)$ tends to $0$ in the size of $G$
assuming that $G$ is non-abelian simple. In this paper we will fix the
finite group $G$ and let $w$ run through $F_{\infty }$.

Alon Amit \cite{amit} has shown that if $G$ is nilpotent then there exists a
constant $c>0$ depending on $G$ only such that for all $w\in $ $F_{\infty }$
we have $P(G,w)>c$. Note that this answers Question 1 affirmatively for
nilpotent groups. He conjectures that the same holds if $G$ is solvable and
that if $G$ is nilpotent then actually 
\begin{equation*}
P(G,w)\geq \frac{1}{\left| G\right| }\text{ for all }w\in F_{\infty }\text{.}
\end{equation*}
He also asked if in turn for a non-solvable finite group $G$ the probability 
$P(G,w)$ can be made arbitrarily small with a suitable $w\in F_{\infty }$.

It is easy to see that Theorem \ref{nagy} already answers Amit's question
affirmatively, but the following stronger result also holds. A group $G$ is 
\emph{just non-solvable} if every proper quotient of $G$ is solvable, but $G$
itself is not.

\begin{theorem}
\label{suru}Let $G$ be a finite just non-solvable group. Then the set 
\begin{equation*}
\left\{ P(G,w)\mid w\in F_{\infty }\right\}
\end{equation*}
is dense in $[0,1]$.
\end{theorem}

\bigskip

\noindent \textbf{Acknowledgement.} The author is grateful to Alon Amit for
communicating his results and questions to him and to Laci Pyber for helpful
advices on how to present the paper.

\section{Proofs}

Let us introduce some notation. Let $G$ be a just non-solvable group and let 
$N$ be a minimal normal subgroup of $G$. Then $N\cong S^{m}$ for some simple
group $S$. By the minimality of $G$ the quotient $G/N$ is solvable, so $S$
is non-abelian and $Z(N)=1$, which implies that $G/Z_{G}(N)$ is non-solvable
so $Z_{G}(N)=1$. Then $G$ embeds into the wreath product 
\begin{equation*}
\mathrm{Aut}(N)\cong \mathrm{Aut}(S)\,\mathrm{wr}\,\mathrm{Sym}(m)
\end{equation*}
where $\mathrm{Sym}(m)$ denotes the symmetric group on $m$ letters and by
the minimality of $N$, $G$ has a transitive image in $\mathrm{Sym}(m)$.
Since $G/N$ is solvable, $N$ is a characteristic subgroup of $G$. Also,
every nontrivial normal subgroup $K\vartriangleleft G$ contains $N$ (using
the minimality of $N$ and that $G$ is just non-solvable). Finally, if in
addition $K=K^{\prime }$ (the commutator subgroup of $K$), then $K=N$.

From now on $G$, $N$, $S$ and $m$ will be as above. Let $G_{j}\cong G$ ($%
1\leq j\leq n$), let 
\begin{equation*}
P=G_{1}\times \cdots \times G_{n}
\end{equation*}
and let 
\begin{equation*}
\pi _{j}:P\rightarrow G_{j}\text{ \ }(1\leq j\leq n)
\end{equation*}
denote the projection to the $j$-th coordinate. Let 
\begin{equation*}
N\cong N_{j}\vartriangleleft G_{j}(1\leq j\leq n)\text{,}
\end{equation*}
let $N_{j}=S_{j,1}\times \cdots \times S_{j,m}$ where $S_{j,i}\cong S$ and
let 
\begin{equation*}
M=N_{1}\times \cdots \times N_{n}\vartriangleleft P\text{.}
\end{equation*}

The first lemma is folklore.

\begin{lemma}
\label{trivi}If $S_{1},\ldots ,S_{n}$ are nonabelian finite simple groups,
then every normal subgroup $K\vartriangleleft S_{1}\times \cdots \times
S_{n} $ is of the form 
\begin{equation*}
K=K_{1}\times \cdots \times K_{n}
\end{equation*}
where $K_{i}=S_{i}$ or $1$ ($1\leq i\leq n$).
\end{lemma}

The next lemma tells us about the normal subgroup structure of subgroups of $%
P$ which project onto each $G_{j}$.

\begin{lemma}
\label{subdir}Let $H\leq P$ be a subgroup containing $M$ such that 
\begin{equation*}
\pi _{j}(H)=G_{j}\text{ \ }(1\leq j\leq n)
\end{equation*}
Let $K$ be a normal subgroup of $H$. Then 
\begin{equation*}
K\cap M=\bigotimes_{\pi _{j}(K)\neq 1}N_{j}
\end{equation*}
\end{lemma}

\begin{proof}
$K\cap M$ is normal in 
\begin{equation*}
M\cong \bigotimes_{1\leq j\leq n\text{, }1\leq i\leq m}S_{j,i}
\end{equation*}
so by Lemma \ref{trivi} it is the direct product of some of the $S_{j,i}$,
that is, 
\begin{equation*}
K\cap M=K_{1}\times \cdots \times K_{n}
\end{equation*}
where $K_{j}\vartriangleleft N_{j}$ ($1\leq j\leq n$).

If $\pi _{j}(K)=1$ then $K\cap N_{j}=1$ so $K_{j}=1$.

If $\pi _{j}(K)\neq 1$ then $\pi _{j}(K)\vartriangleleft \pi _{j}(H)=G_{j}$
so $N_{j}\subseteq \pi _{j}(K)$, since $N_{j}$ is a minimal normal subgroup
in $G_{j}$. In this case 
\begin{equation*}
K\cap M\supseteq \lbrack K,M]\supseteq \lbrack K,N_{j}]=[\pi
_{j}(K),N_{j}]=N_{j}
\end{equation*}
so $K_{j}=N_{j}$ (here we use the direct product form and that the
commutator $[N_{j},N_{j}]=N_{j}$).

The lemma holds.
\end{proof}

\bigskip

Let 
\begin{equation*}
\left( a_{1},\ldots ,a_{k}\right) ,\text{ }\left( b_{1},\ldots ,b_{k}\right)
\in G^{k}
\end{equation*}
be $k$-tuples from $G$. We say that $\left( a_{1},\ldots ,a_{k}\right) $ and 
$\left( b_{1},\ldots ,b_{k}\right) $ are \emph{automorphism independent over 
}$G$ if there exists no $\alpha \in \mathrm{Aut}(G)$ such that $%
a_{i}^{\alpha }=b_{i}$ for all $1\leq i\leq k$.

Our next lemma shows that subgroups of $G_{1}\times \cdots \times G_{n}$
satisfying some natural conditions contain $N_{1}\times \cdots \times N_{n}$.

\begin{lemma}
\label{trukk}Let $a_{i,j}\in G_{j}$ ($1\leq i\leq k$, $1\leq j\leq n$), such
that we have 
\begin{equation*}
\left\langle a_{1,j},\ldots ,a_{k,j}\right\rangle =G_{j}\text{ }(1\leq j\leq
n)
\end{equation*}
and that for all $1\leq j<l\leq n$ the $k$-tuples $\left( a_{1,j},\ldots
,a_{k,j}\right) $ and $\left( a_{1,l},\ldots ,a_{k,l}\right) $ are
automorphism independent over $G$. For $1\leq i\leq k$ let 
\begin{equation*}
h_{i}=\left( a_{i,1},\ldots ,a_{i,n}\right) \in G_{1}\times \cdots \times
G_{n}
\end{equation*}
and let 
\begin{equation*}
H=\left\langle h_{1},\ldots ,h_{k}\right\rangle \leq G_{1}\times \cdots
\times G_{n}
\end{equation*}
Then 
\begin{equation*}
M=N_{1}\times \cdots \times N_{n}\leq H.
\end{equation*}
\end{lemma}

\begin{proof}
Let 
\begin{equation*}
f:G_{1}\times \cdots \times G_{n}\rightarrow G_{1}\times \cdots \times
G_{n-1}
\end{equation*}
denote the projection to the first $n-1$ coordinates. Let $H_{1}=f(H)$ and
let 
\begin{equation*}
R=\pi _{n}(\mathrm{Ker}(f))\leq G_{n}\text{.}
\end{equation*}
By induction on $n$, we have $N_{1}\times \cdots \times N_{n-1}\leq H_{1}$.
Also $R$ is normal in $G_{n}$ so by the minimality of $N_{n}$ in $G_{n}$
either $N_{n}\leq R$ or $R=1$.

We claim that $N_{n}\leq R$. Assume $R=1$. Let us define the function $%
\varphi :H_{1}\rightarrow G_{n}$ by 
\begin{equation*}
\varphi (g_{1},\ldots ,g_{n-1})=g_{n}\text{ if }(g_{1},\ldots
,g_{n-1,}g_{n})\in H\text{. }
\end{equation*}
Then $\varphi $ is well-defined, since 
\begin{equation*}
(g_{1},\ldots ,g_{n-1,}g_{n}),(g_{1},\ldots ,g_{n-1,}g_{n}^{\prime })\in H
\end{equation*}
implies $g_{n}^{-1}g_{n}^{\prime }\in R$. So $\varphi $ is a homomorphism.
Using $\left\langle a_{1,n},\ldots ,a_{k,n}\right\rangle =G_{n}$ we also see
that $\varphi $ is surjective.

Let $K=\mathrm{Ker}(\varphi )$. Then $K$ is normal in $H_{1}$ and 
\begin{equation*}
H_{1}/K\cong G_{n}\cong G
\end{equation*}
which is not solvable. Since $N_{1}\times \cdots \times N_{n-1}\subseteq
H_{1}$, the use of Lemma \ref{subdir} for $H_{1}$ and $K$ gives us 
\begin{equation*}
K\cap M=\bigotimes_{\pi _{j}(K)\neq 1}N_{j}\text{. }
\end{equation*}

Now $M\subseteq K$ would imply that $H_{1}/K$ is solvable, a contradiction.
So there exists a coordinate $1\leq l<n$ such that $\pi _{l}(K)=1$, that is, 
$K\subseteq \mathrm{Ker}(\pi _{l})$. Moreover 
\begin{equation*}
H_{1}/\mathrm{Ker}(\pi _{l})\cong G_{l}\cong G
\end{equation*}
which implies $K=\mathrm{Ker}(\pi _{l})$. This shows that the function $%
\alpha :G_{l}\rightarrow G_{n}$ defined by 
\begin{equation*}
\alpha (g_{l})=g_{n}\text{ if }(g_{1},\ldots ,g_{l},\ldots ,g_{n})\in H
\end{equation*}
is an isomorphism. In particular, $\alpha (a_{i,l})=a_{i,n}$ ($1\leq i\leq k$%
), so the $k$-tuples $\left( a_{1,l},\ldots ,a_{k,l}\right) $ and $\left(
a_{1,n},\ldots ,a_{k,n}\right) $ are not automorphism independent over $G$
which contradicts the assumptions of the lemma. So the claim $N_{n}\leq R$
holds and so $1\times \cdots \times 1\times N_{n}\leq H$.

Now let $L=f^{-1}(N_{1}\times \cdots \times N_{n-1})\leq H$. Let $L^{(i)}$
denote the $i$-th element of the derived series of $L$ and let $r$ be a
number such that $L^{(r)}=L^{(r+1)}$. Then $f(L^{(r)})=N_{1}\times \cdots
\times N_{n-1}$ and since $1\times \cdots \times 1\times N_{n}\leq L$ also $%
1\times \cdots \times 1\times N_{n}\leq L^{(r)}$. Now $J=\pi _{n}(L^{(r)})$
is normal in $G_{n}$, $N_{n}\leq J$ and $J^{\prime }=J$, so $J=N_{n}$. This
implies 
\begin{equation*}
L^{(r)}=N_{1}\times \cdots \times N_{n}\leq H
\end{equation*}
what we wanted to prove.
\end{proof}

\bigskip

\noindent \textbf{Remark.} This lemma is well-known in the subcase when $G$
is a nonabelian finite simple group (see \cite{canlub}, or \cite{wie}). We
will state a light corollary of that which we will use in the proof of
Theorem \ref{nagy}.

\begin{corollary}
\label{easy}Let $G_{i}$ ($1\leq i\leq n$) be finite nonabelian simple groups
and let 
\begin{equation*}
H\leq G_{1}\times \cdots \times G_{n}
\end{equation*}
such that the projections $\pi _{i}(H)=G_{i}$ ($1\leq i\leq n$). Then there
exists $g\in H$ such that $\pi _{i}(g)\neq 1$ ($1\leq i\leq n$).
\end{corollary}

\begin{proof}
We proceed by induction on $n$. For $n=1$ the lemma is trivial. By induction
we have an element $g\in H$ such that $\pi _{i}(g)\neq 1$ ($1\leq i<n$). If
the last coordinate is automorphism dependent on some previous coordinate $k$
then $\pi _{k}(g)\neq 1$ implies $\pi _{n}(g)\neq 1$. If it is not, then $%
1\times \cdots \times 1\times G_{n}\leq H$ and we can set the last
coordinate of $g$ as we wish.
\end{proof}

\bigskip

Let $S$\ be a set of finite simple groups. We say that a finite group $G$ is
in $\mathrm{Comp}(S)$ if all nonabelian composition factors of $G$ are in $S$%
. An affirmative answer for the following question would be a far-reaching
generalization of Theorem \ref{nagy}.

\bigskip

\noindent \textbf{Question 2. }\textit{Let }$S$\textit{\ be a finite set of
finite simple groups. Is it true that for all }$n$\textit{\ there exists }$%
w\in F_{n}$\textit{\ such that every quotient of the one-relator group }$%
F_{n}/\left\langle w\right\rangle $\textit{\ which lies in }$\mathrm{Comp}%
(S) $\textit{\ is solvable?}

\bigskip

Note that we dont know the answer even in the case when $S$ consists of one
simple group.

\bigskip

Now we prove Theorem \ref{suru}.

\bigskip

\begin{proof}[Proof of Theorem 3]
Let $m$ be the number of maximal subgroups of $G$. Let $d>\log _{2}m$ be an
integer to be chosen later. The probability that $d$ independent random
elements all fall into a fixed maximal subgroup $M$ is at most $\left|
G:M\right| ^{-d}\leq 2^{-d}$ so the probability that $d$ random elements do
not generate $G$ is at most $m2^{-d}<1$. In particular, $G$ can be generated
by $d$ elements. Let 
\begin{equation*}
S=\left\{ (g_{1},\ldots ,g_{d})\in G^{d}\mid g_{1},\ldots ,g_{d}\text{
generate }G\right\}
\end{equation*}
be the set of generating $d$-tuples.

Now $\mathrm{Aut}(G)$ acts on $S$ by $(g_{1},\ldots ,g_{d})^{\alpha
}=(g_{1}^{\alpha },\ldots ,g_{d}^{\alpha })$ where $\alpha \in \mathrm{Aut}%
(G)$. This action is fixed-point free, as if $\alpha $ fixes all the
elements of a generating set then it fixes the whole $G$. Let $r$ be the
number of $\mathrm{Aut}(G)$-orbits and let $t_{1},\ldots ,t_{r}\in S$ be an
orbit representative system.

It is easy to see that the conditions of Lemma \ref{trukk} hold for $%
a_{i,j}=t_{j}(i)$. This implies that the $r$-tuples 
\begin{equation*}
h_{i}=\left( t_{1}(i),\ldots ,t_{n}(i)\right)
\end{equation*}
generate a group $H$ which contains $N_{1}\times \cdots \times N_{r}$.

Let $1\neq g\in N$, let $k\leq r$ be a natural number to be chosen later and
let the $r$-tuple $h$ be defined by 
\begin{equation*}
h(i)=1\text{ (}1\leq i\leq k\text{) and }h(i)=g\text{ (}k<i\leq r\text{)}
\end{equation*}
Then $h\in N_{1}\times \cdots \times N_{r}\subseteq H$, so there exists a
word $w\in F_{d}$ such that $w(h_{1},\ldots ,h_{d})=h$.

Now let us evaluate $w$ on the set of possible $d$-tuples from $G$. We
completely control the evaluation on generating tuples; since 
\begin{equation*}
w(g_{1}^{\alpha },\ldots ,g_{d}^{\alpha })=(w(g_{1},\ldots ,g_{d}))^{\alpha }%
\text{ (}\alpha \in \mathrm{Aut}(G)\text{)}
\end{equation*}
we have 
\begin{equation*}
\left| \left\{ (g_{1},\ldots ,g_{d})\in S\mid w(g_{1},\ldots
,g_{d})=1\right\} \right| =k\left| \mathrm{Aut}(G)\right| \text{.}
\end{equation*}

On $d$-tuples $(g_{1},\ldots ,g_{d})$ not generating $G$ we do not control $%
w(g_{1},\ldots ,g_{d})$. This gives the estimate 
\begin{equation*}
k\frac{\left| \mathrm{Aut}(G)\right| }{\left| G\right| ^{d}}\leq P(G,w)\leq k%
\frac{\left| \mathrm{Aut}(G)\right| }{\left| G\right| ^{d}}+m2^{-d}
\end{equation*}
and for $k=r$ we get 
\begin{equation*}
k\frac{\left| \mathrm{Aut}(G)\right| }{\left| G\right| ^{d}}=\frac{\left|
S\right| }{\left| G\right| ^{d}}\geq 1-m2^{-d}\text{.}
\end{equation*}

Since $d$ can be chosen to be arbitrarily large, both $\left| \mathrm{Aut}%
(G)\right| /\left| G\right| ^{d}$ and $m2^{-d}$ get arbitrarily small. Now $%
k\leq r$ is arbitrary which shows that the set 
\begin{equation*}
\left\{ P(G,w)\mid w\in F_{\infty }\right\}
\end{equation*}
is dense in $[0,1]$.
\end{proof}

\bigskip

The answer to Amit's question follows as an easy corollary of Theorem \ref
{suru}.

\begin{corollary}
Let $G$ be a finite non-solvable group. Then the set 
\begin{equation*}
\left\{ P(G,w)\mid w\in F_{\infty }\right\}
\end{equation*}
accumulates in $0$.
\end{corollary}

\begin{proof}
Let $K$ be a normal subgroup in $G$ such that $G/K$ is just non-solvable and
let $g_{1},\ldots ,g_{n}$ be independent uniform random elements of $G$.
Then $g_{1}K,\ldots ,g_{n}K$ are independent uniform random elements of $G/K$
which yields 
\begin{equation*}
P(G/K,w)=P(w(g_{1},\ldots ,g_{n})\in K)\geq P(w(g_{1},\ldots
,g_{n})=1)=P(G,w)
\end{equation*}
for $w\in F_{\infty }$. Using Theorem \ref{suru} we get that for every $%
\epsilon >0$ we have $w\in F_{\infty }$ such that 
\begin{equation*}
P(G,w)\leq P(G/K,w)<\epsilon
\end{equation*}
and so the corollary holds.
\end{proof}

\bigskip

We are ready to prove Theorem \ref{nagy}.

\bigskip

\begin{proof}[Proof of Theorem 1]
For each subgroup $H\leq G$ let us choose a homomorphism 
\begin{equation*}
\varphi _{H}:H\rightarrow .
\end{equation*}
as follows. If $H$ is solvable then let $\varphi _{H}=\mathrm{Id}$ be the
identity, otherwise let $\varphi _{H}$ be a homomorphism to a just
non-solvable quotient of $H$.

Let us enumerate all the $n$-tuples from $G$ as $t_{1},t_{2},\ldots ,t_{k}$
where $k=\left| G\right| ^{n}$. Let $t_{i,j}$ denote the $j$-th element of $%
t_{i}$ ($1\leq j\leq n$). For $1\leq i\leq k$ let 
\begin{equation*}
H_{i}=\left\langle t_{i,1},t_{i,2},\ldots ,t_{i,n}\right\rangle
\end{equation*}
Let $\varphi _{i}=\varphi _{H_{i}}$ and let $G_{i}=\varphi _{i}(H_{i})$. Let 
$N_{i}$ be the minimal normal subgroup of $G_{i}$ if $G_{i}$ is just
non-solvable, otherwise let $N_{i}=1$. Also let 
\begin{equation*}
u_{i,j}=\varphi _{i}(t_{i,j})\text{ }(1\leq i\leq k,1\leq j\leq n)
\end{equation*}
and let 
\begin{equation*}
p_{j}=(u_{1,j},u_{2,j},\ldots ,u_{k,j})\in G_{1}\times \cdots \times G_{k}%
\text{ }(1\leq j\leq n)
\end{equation*}
Let 
\begin{equation*}
L=\left\langle p_{1},p_{2},\ldots ,p_{n}\right\rangle \leq G_{1}\times
\cdots \times G_{k}
\end{equation*}
and let 
\begin{equation*}
\pi _{i}:G_{1}\times \cdots \times G_{k}\rightarrow G_{i}\text{ }(1\leq
i\leq k)
\end{equation*}
denote the projection to the $i$-th coordinate. Then $\pi _{i}(L)=G_{i}$ ($%
1\leq i\leq k$). Let $L^{(i)}$ denote the $i$-th derived subgroup of $L$ and
let $r$ be an integer such that $M=L^{(r)}=L^{(r+1)}$. Then $\pi
_{i}(M)\vartriangleleft G_{i}$ and $\pi _{i}(M)^{\prime }=\pi _{i}(M)$ so $%
\pi _{i}(M)=N_{i}$. Now all the $N_{i}\neq 1$ are isomorphic to some direct
power of a nonabelian simple group so $M$ lies in a direct product of
nonabelian simple groups and projects to each factor of the product. By
Corollary \ref{easy} there exists an element $g\in M\leq L$ such that $\pi
_{i}(g)\neq 1$ if and only if $N_{i}\neq 1$. Let $w\in F_{n}$ be a word such
that $w(p_{1},p_{2},\ldots ,p_{n})=g$.

We claim that this $w$ will be good for our purposes. Indeed, we have 
\begin{equation*}
\pi _{i}(g)=w(u_{i,1},\ldots ,u_{i,n})=w(\varphi _{i}(t_{i,1}),\ldots
,\varphi _{i}(t_{i,n}))=\varphi _{i}(w(t_{i,1},\ldots ,t_{i,n}))
\end{equation*}
Now if $H_{i}$ is solvable then $\varphi _{i}$ is the identity map and $\pi
_{i}(g)=1$, so we get $w(t_{i,1},\ldots ,t_{i,n})=1$. If $H_{i}$ is not
solvable, then $\pi _{i}(g)\neq 1$ and since $\varphi _{i}$ is a
homomorphism we have $w(t_{i,1},\ldots ,t_{i,n})\neq 1$. The theorem holds.
\end{proof}

\bigskip

\noindent \textbf{Remark on Question 1.} Let $G$ be a finite group for which
there exists a constant $c>0$ such that for all $w\in $ $F_{\infty }$ we
have $P(G,w)>c$. Then as we saw for large enough $d$ most of the $d$-tuples
generate $G$ and so for every word $w\in F_{d}$ there exists a generating
set $\left\langle g_{1},\ldots ,g_{d}\right\rangle =G$ such that $%
w(g_{1},\ldots ,g_{d})=1$, that is, $G$ is a quotient of the one-relator
group $F_{d}/\left\langle w\right\rangle $. In particular, Amit's result 
\cite{amit} implies an affirmative answer for Question 1 for finite
nilpotent groups.

\end{document}